\newcommand{\be}{\begin{equation}} \newcommand{\ee}{\end{equation}}
\newcommand{\bea}{\begin{eqnarray}}\newcommand{\eea}{\end{eqnarray}}
\newcommand{\nn}{\nonumber}
\newcommand{\ba}{\begin{array}}\newcommand{\ea}{\end{array}}
\newcommand{\R}{\mbox{\bf R}}\newcommand{\U}{\mbox{\bf U}}

\newcommand\Tq{\mbox{Tr}_q}

\documentstyle[12pt]{article}
\begin{document}
\renewcommand{\thefootnote}{\fnsymbol{footnote}}
\begin{center}
{\large QUANTUM DEFORMATIONS OF MULTI-INSTANTON SOLUTIONS  \\
IN THE TWISTOR SPACE\\
{\it B.M.Zupnik }}\\
{\it Bogoliubov   Laboratory of Theoretical Physics , JINR, Dubna,
Moscow Region, 141980, E-mail: zupnik@thsun1.jinr.dubna.su }\\

\vspace{0.5cm}
To be published in "Pis'ma v ZhETP" v.62, n.4
\end{center}

We consider the quantum-group self-duality equation in the framework
 of the gauge theory on a deformed twistor space. Quantum deformations
of the Atiyah-Drinfel'd-Hitchin-Manin and t'Hooft multi-instanton
solutions are constructed.

\vspace{0.5cm}
\renewcommand{\thefootnote}{\arabic{footnote}}
\setcounter{footnote}{0}

   The quantum-group gauge theory was considered in the framework of
 the algebra of local differential complexes \cite{a1}-\cite{a3}
 or as a noncommutative generalization of the fibre bundles over the
classical or quantum basic spaces \cite{a4,a5}.

   We prefer to use  local constructions of the noncommutative
connection forms or gauge fields as a deformed analogue of the local
gauge  fields. In particular, the quantum-group self-duality equation
(QGSDE)
has been considered in the deformed 4-dimensional Euclidean space, and
an explicit formula for the corresponding one-instanton solution has
been constructed \cite{a3}. This
solution can be treated as $q$-deformation of the BPST-instanton
\cite{a6}. We shall discuss here  quantum deformations of the general
multi-instanton solutions \cite{a7}.

The conformal covariant description of the classical ADHM solution was
considered in Ref\cite{a8}. We shall study the quantum deformation of
this version of the twistor formalism. It is convenient to
discuss firstly the deformations of the complex conformal group
$GL(4,C)$ , complex twistors and the complex linear gauge groups.

Let $R^{ab}_{cd},\; (a, b, c, d\ldots =1\ldots 4)$ be the solution of
the 4D Yang-Baxter equation satisfying also the Hecke relation
\bea
& R\;R^\prime \;R = R^\prime \;R\;R^\prime & \\ \label{A1}
& R^2 = I + (q-q^{-1})R  & \label{A2}
\eea
where $q$ is a complex parameter. Note that the standard notation for
these $R$-matrices is $ R=\hat{R}_{12},\;R^\prime=\hat{R}_{23}$
\cite{a9}.

Consider also the $SL_q(2,C)\;R$-matrix
\be
R^{\alpha\beta}_{\mu\nu}=q\delta^\alpha_\mu \delta^\beta_\nu +
\varepsilon^{\alpha\beta}(q)\varepsilon_{\mu\nu}(q) \label{A3}
\ee
where $\varepsilon(q)$ is the deformed antisymmetric symbol.

Noncommutative twistors were considered in Ref\cite{a10}. We shall
use the $R$-matrix approach to define the differential calculus on
the deformed twistor space.

Let $z^\alpha_a$ and $dz^\alpha_a$ be the components of the $q$-twistor
and their differentials
\bea
& R^{\alpha\beta}_{\mu\nu}\;z^\mu_a\; z^\nu_b =z^\alpha_c\; z^\beta_d\;
R^{dc}_{ba} &\\ \label{A4}
& z^\alpha_a\;dz^\beta_b = R^{\alpha\beta}_{\mu\nu}\;dz^\mu_c\; z^\nu_d\;
R^{dc}_{ba} &\\ \label{A5}
& dz^\alpha_a\;dz^\beta_b = -R^{\alpha\beta}_{\mu\nu}\;dz^\mu_c\;dz^\nu_d
\;R^{dc}_{ba} & \label{A6}
\eea

One can define also the algebra of partial derivatives
$\partial^a_\alpha$
\bea
& R^{ab}_{cd}\;\partial^c_\alpha\; \partial^d_\beta=\partial^a_\mu\;
\partial^b_\nu\; R_{\beta\alpha}^{\nu\mu} & \\ \label{A7}
& \partial^a_\alpha\; z^\beta_b =\delta^a_b\; \delta_\alpha^\beta +
R^{\beta\mu}_{\alpha\nu}\; R^{da}_{cb}\; z^\nu_d\; \partial^c_\mu &
\label{A8}
\eea

Consider the 4D deformed $\varepsilon_q$-symbol
\be
R^{ba}_{fe}\; \varepsilon_q^{efcd}=-\frac{1}{q}
\varepsilon_q^{abcd} \label{A9}
\ee

The $q$-twistors satisfy the following identity:
\be
\varepsilon_q^{abcd}\; z^\beta_b\; z^\mu_c\; z^\nu_d = 0 \label{A10}
\ee

The $SL_q(2)$-invariant bilinear function of twistors has the zero length
in the projective 6D vector space
\bea
& y_{ab}=\varepsilon_{\alpha\beta}(q)\; z^\alpha_a\; z^\beta_b =
[P^{(-)}]^{dc}_{ba}\;y_{cd} & \\ \label{A11}
& (y,y)=\varepsilon_q^{abcd}\; y_{ab}\; y_{cd}=0 & \label{A12}
\eea

Consider a duality transformation $\ast$ of the basic differential
2-forms
 \cite{a3}
\be
\ast dz\; dz^\prime = dz\; dz^\prime\; P^{(+)} - dz\; dz^\prime\; P^{(-)}
\label{A13}
\ee
where $P^{(\pm)}$ are the projection operators of $GL_q(4)\;$\cite{a9}.
Note that the self-dual part $dz\; dz^\prime\; P^{(+)}$ is proportional
to
\be
\varepsilon_{\alpha\beta}(q)\; dz^\alpha_a\; dz^\beta_b \label{A14}
\ee

Let $T^{i}_{k}$  be matrix elements of the
$GL_q(N)$ quantum group
\be
R_G T T^\prime =T T^\prime R_G \label{A15}
\ee
where $R_G$ is the $R$-matrix of $GL_q(N)$.

Quantum deformation of the $GL_q(N)$ gauge connection can be treated
in terms of the noncommutative algebra for the components $A^{i}_{k}$
of the connection 1-form \cite{a1,a2}
\be
(A\; R_G\; A + R_G\; A\; R_G\; A\; R_G )^{ikl}_{mnp}=0 \label{A16}
\ee
where $i, k, l, m, n, p=1\ldots N$.
These relations generalize the anticommutativity conditions for
components of the classical connection form.

The restriction on the quantum trace of the connection $\alpha= \Tq A=0$
is inconsistent with (\ref{A16}),
 but we can use the gauge-covariant relations
$\alpha^2=0,\;\Tq A^2=0\;\mbox{and}\;d\alpha=0\; $\cite{a3}.
 The curvature 2-form
$F = dA - A^2$ is $ q$-traceless for this model.

Consider the explicit realization of this gauge algebra in terms of
$z, dz $ and the set $B$ of additional noncommutative parameters
\be
A^{i}_{k}(z,dz,B)=dz^\alpha_a\; A^{ai}_{\alpha k}(z,B) \label{A17}
\ee
The analogous realizations were considered on the $GL_q(2)$ and $E_q(4)$
quantum spaces \cite{a1}-\cite{a3}. We shall treat the representation
 (\ref{A17}) as a local gauge field on the $q$-twistor space.

Let us consider the quantum deformation of the $GL(2)$ t'Hooft solution
\cite{a8}
\bea
&A^\alpha_\beta =q^{-3}dz^\alpha_a\; (\partial^a_\mu \Phi)\Phi^{-1}
\varepsilon^{\sigma\mu}(q)\varepsilon_{\sigma\beta}(q)&\\ \label{A18}
& \Phi =\sum_{i} (X^i)^{-1},\;\;\;\;\;X^i =(y,b^i)= \varepsilon_q^{abcd}
 y_{ab}\;b^i_{cd} & \label{A19}
\eea
where $b^i_{cd}$ are the noncommutative isotropic $6D$ vectors
\bea
& db^i_{cd}=0,\;\;\;\;\;(b^i,b^i)=0&\\ \label{A20}
& [y_{ab}, X^i ]=[b^i_{cd},X^i ]=0 & \label{A21}
\eea

The central elements $X^i$ of the (B,z)-algebra do not commute with $dz$
\be
X^i\; dz^\alpha_a = q^2 dz^\alpha_a\; X^i \label{A22}
\ee

Stress that $A^\alpha_\beta$ satisfies  Eq(\ref{A16}) and its quantum
trace is a $U(1)$-gauge field with the zero field-strength
\be
\Tq A =-q^{-3}d\Phi \Phi^{-1},\;\;\;\;\; \Tq dA=0 \label{A23}
\ee

The QGSDE for $A^\alpha_\beta$ is equivalent to the finite-difference
Laplace equation for the function $\Phi$ on the $q$-twistor space
\bea
& \Delta^{ba} \Phi(X^i)=\sum_{i}\Delta^{ba}\frac{1}{X^i} =0&\\\label{A24}
& \Delta^{ba}\Phi=\frac{q}{1+q^2}\varepsilon^{\alpha\beta}(q)
\partial^b_\beta
\partial^a_\alpha \Phi=(\partial^{ba} +\frac{1}{2} y_{cd}\partial^{dc}
\partial^{ba})\Phi&\\
\label{A25}
&\partial^{ba} y_{cd}=[P^{(-)}]^{ab}_{dc},\;\;\;\;\;
\partial^{ba}(X^i)^{-1}=-q^{-2}(X^i)^{-2}(b^i )^{ab}& \label{A26}
\eea

 The ADHM-twistor functions of Ref\cite{a7} can be connected with some
$GL(N+2k)$ matrix function. Let us introduce the notation for indices
of different types: $I, K, L, M= 1\ldots N+2k$ and $A, B =1\ldots k$.
The Ansatz for the general self-dual $GL_q(N,C)$ field contains the
deformed twistors $u(z)$ and $\tilde{u}(z)$
\be
A^i_k =du^i_I \tilde{u}^I_k,\;\;\;\; u^i_I \tilde{u}^I_k=\delta^i_k
\label{A27}
\ee

The commutation relations for the $u$ and $ \tilde{u}$ twistors are
\bea
& (R_G)^{ik}_{lm}u^l_I u^m_K = u^i_L u^k_M \R^{LM}_{IK} &\\\label{A28}
& \R^{KI}_{ML}\tilde{u}^L_i \tilde{u}^M_k =\tilde{u}^I_l \tilde{u}^K_m
 (R_G)_{ki}^{ml} & \\ \label{A29}
& \tilde{u}^I_l (R_G)^{li}_{mk} u^m_K = u^i_L \R^{IL}_{KM}\tilde{u}^M_k&
 \label{A30}
\eea
where the $R$-matrices for $GL_q(N,C)$ and $GL_q(N+2k,C)$ are used.

Consider also the linear twistor functions $v$ and $\tilde{v}$
\bea
& v^{A\alpha}_I = z^\alpha_a b^{aA}_I &\\ \label{A31}
& \tilde{v}^{IA\alpha} =\tilde{b}^{IAa} z^\alpha_a & \label{A32}
\eea

Introduce the following condition for these functions:
\be
v^{A\alpha}_I \tilde{v}^{IB\beta}=g^{AB}(z)\varepsilon^{\alpha\beta}(q)
\label{A33}
\ee
where $g(z)$ is the nondegenerate $(k\times k)$ matrix with the central
elements
\be
g^{AB}(z)=\frac{q}{1+q^2}b^{aA}_I\; \tilde{b}^{IBb}\; y_{ab} \label{A34}
\ee

The condition (\ref{A33}) is equivalent to the restriction on the
elements of the $B$-algebra
\be
[P^{(+)}]_{ab}^{cd}\; b^{aA}_I\; \tilde{b}^{IBb} =0 \label{A35}
\ee

Write the basic commutation relations of the $B$-algebra
\bea
& R^{ab}_{cd}\; b^{cA}_I\; b^{dB}_K =b^{aB}_L\; b^{bA}_M\;\R^{ML}_{KI}&\\
\label{A36}
& \R_{LM}^{IK}\;\tilde{b}^{LAa}\;\tilde{b}^{MBb}=R^{ab}_{cd}\;
\tilde{b}^{IBc}\;
\tilde{b}^{KAd} &\\ \label{A37}
& R^{ab}_{cd}\;b^{cA}_I\; \tilde{b}^{KBd}=\R^{KL}_{IM}\;\tilde{b}^{MBa}\;
b^{aB}_L
& \label{A38}
\eea
Remark that a formal permutation of the indices $A$ and $B$ is
commutative. It is not difficult to define the relations between $b,
\tilde{b}$ and $z, dz$.

Consider the new functions
\be
\tilde{v}_{A\alpha}^I=g_{AB}(z)\varepsilon_{\alpha\beta}(q)
\tilde{v}^{IB\beta} \label{A39}
\ee
where we use the inverse matrix with respect to the matrix (\ref{A34}).

Now one can construct the full quantum $GL_q(N+2k,C)$ matrices
\be
\U=\left( u^i_I\atop v^{A\alpha}_I \right)\;,\;\;\;\;\;\; \U^{-1}=\left(
\tilde{u}^I_i \atop\tilde{v}_{A\alpha}^I \right)
\label{A40}
\ee
The standard $GL_q(N+2k,C)$ commutation relations for these matrices
contain Eqs(28-30) and the relations for the $v$ and $\tilde{v}$
functions.

Write explicitly the orthogonality and completeness conditions for the
deformed ADHM-twistors
\bea
& u^i_I\;\tilde{v}^{IA\alpha} =0 &\\ \label{A41}
& v^{A\alpha}_I \;\tilde{u}^I_i=0 &\\ \label{A42}
& \delta^I_K =\tilde{u}^I_i\;u^i_K\;+\;\tilde{v}^{IA\alpha}\;
g_{AB}(z)\varepsilon_{\alpha\beta}(q)\;v^{B\beta}_K &\label{A43}
\eea

The gauge-field algebra (\ref{A16}) for the deformed ADHM-Ansatz
(\ref{A27}) can be generated by the differential algebra on the
$GL_q(N+2k,C)$ matrices $\U,\;\U^{-1},\;d\U$ which contains the following
relations:
\bea
& \tilde{u}^I_i\;(R_G )^{ik}_{lm}\;du^l_K = du^k_L\;(\R^{-1})^{IL}_{KM}
\tilde{u}^M_m &\\ \label{A44}
& du^i_L\;du^k_M\;(\R^{-1})^{LM}_{IK}= -(R^{-1}_G )^{ik}_{lm}\;
du^l_I\;du^m_K & \label{A45}
\eea
These relations are consistent with the commutation relation (28-30).

The self-duality of the connection (27) follows from  Eqs(31,32,41-43).
\bea
& dA^i_k - A^i_l\;A^l_k = du^i_I\;(\tilde{u}^I_l\;u^l_M\;-\;\delta^I_M)
d\tilde{u}^M_k = &\\ \nn
& =-u^i_I\;\tilde{b}^{IAa}\;g_{AB}(z)\varepsilon_{\alpha\beta}(q)\;
dz^\alpha_a\;dz^\beta_b\;b^{Bb}_M\;\tilde{u}^M_k & \label{A46}
\eea
This curvature contains the self-dual 2-form (14) only.

It should be stressed that all $R$-matrices of our deformation scheme
satisfy the Hecke relation with the common parameter $q$. The other
possible parameters of different $R$-matrices are independent. The case
$q=1$ corresponds to the unitary deformations $(R^2=I) $ of the twistor
space and the gauge groups. It is evident that the trivial deformation
of the $z$-twistors is consistent with the nontrivial unitary deformation
of the gauge sector and vice versa.

The Euclidean conformal $q$-twistors are a representation of the
$U^* (4)\times SU_q (2) $ group. The antiinvolution for these twistors
has the following form:
\be
(z^\alpha_a )^* = \varepsilon_{\alpha\beta}(q)\;z^\beta_b\;C^b_a
\label{A47}
\ee
where $C$ is the charge conjugation matrix for $U^* (4)$. We can use the
gauge group $U_q (N)$ in the framework of our approach.

An analogous construction can be considered for the real twistors and the
gauge group $GL_q (N,R)$.

The author would like to thank
A.T. Filippov, E.A. Ivanov, A.P.Isaev and
V.I. Ogievetsky
 for helpful discussions and interest in this work.

I am grateful to administration of JINR and Laboratory of Theoretical
Physics for hospitality. This work was supported in part by
 ISF-grant RUA000, INTAS-grant 93-127  and the contract No.40 of Uzbek
 Foundation of Fundamental
Research .
\renewcommand{\refname}{}\refname

\end{document}